\documentclass{IEEEtran}
        
\usepackage{graphicx}
\usepackage[dvipsnames]{xcolor}
\usepackage{mathtools}

\title{A Brief Introduction to  Quantum Network Control} %\\  and Open Problems}
\author{
    \IEEEauthorblockN{V\'ictor Valls\IEEEauthorrefmark{1}, Panagiotis Promponas\IEEEauthorrefmark{2}, Leandros Tassiulas\IEEEauthorrefmark{2}} \\
    \IEEEauthorblockA{\IEEEauthorrefmark{1}IBM Research Europe -- Dublin
    }
    \IEEEauthorblockA{\IEEEauthorrefmark{2}Department of Electrical Engineering, Yale University
    \\}
}

\begin{document}
\maketitle

\begin{abstract}

Quantum networking is an emerging area with the potential to transform information processing and communications. In this paper, we present a brief introduction to quantum network control, an area in quantum networking dedicated to designing algorithms for distributing entanglement (i.e., entangled qubits). We start by explaining what qubits and entanglement are and how they furnish quantum network control operations such as entanglement swapping and teleportation. With those operations, we present a model for distributing entanglement in a multi-hop quantum network to enable applications such as quantum key distribution and distributed quantum computing. We conclude the paper by presenting open research problems in the field, including the characterization of the quantum network capacity region and the design of throughput-optimal policies.

\end{abstract}

%%%%%%%%%%%%%%%%%%%%%%%%
\section{Introduction}

Quantum networks are expected to transform information processing and communication by leveraging the unique properties of quantum mechanics. These networks can enhance data security and enable distributed quantum computing applications to solve problems beyond the reach of current computer systems.

However, quantum networks are still in their infancy. Single-hop communications are challenging in practice, and many technical issues must be addressed before quantum network deployments become a reality \cite{castelvecchi2018quantum, WEH18}.  Nonetheless, the building blocks of how quantum networks will operate are starting to take shape, which has motivated researchers to start developing the algorithms to operate them when the hardware is available \cite{Pant:2019vq, 9951209}.

This article presents an introduction to \emph{quantum network control}: an area in quantum networking that focuses on designing algorithms for distributing entanglement (i.e., entangled qubits) across a multi-hop quantum network. Entanglement is key for realizing many quantum networking applications such as Quantum Key Distribution (QKD). This article is specifically aimed at researchers in the field of communications without prior experience of quantum networking, quantum physics, or quantum computing.  One of the challenges when starting with quantum computing/networking is the terminology used in physics. Concepts such as qubit or entanglement are not hard to understand on their own. However, when they are put together to do something  different from what we are used to in classical networking, the full picture becomes less clear. Or, put differently, what do qubits and entanglements allow us to do that we could not do before? Answering this question is crucial for understanding what quantum networks can do, but also what they cannot. A common misconception is that quantum networks will allow us to transmit large volumes of data at a rate that exceeds the speed of light \cite{BBVA}. However, that could not be further from the truth. Quantum networks cannot transmit classical or quantum information faster than the speed of light. Instead, quantum networks will work in tandem with existing communication networks to enable applications that were not possible before, such as QKD, distributed quantum computing, and quantum sensing \cite{castelvecchi2018quantum}.

Understanding what quantum networks can or cannot do is essential for designing algorithms that can effectively operate them. Therefore, we begin this article by explaining a simplified quantum network application, which will highlight the key differences between quantum and traditional communication networks. Next, we explain the fundamental concepts of qubits and entanglement, showing how they enable important quantum control operations such as entanglement swapping and quantum teleportation. Using these operations, we introduce a basic quantum network control model and discuss open research questions in the field, including the characterization of the quantum network capacity region and the design of throughput optimal policies for entanglement distribution.

\begin{figure}
\centering
\includegraphics[width=0.9\columnwidth]{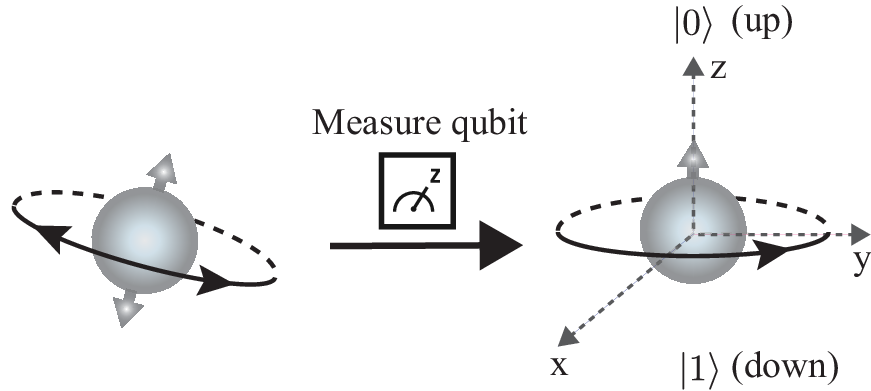}
\caption{A qubit in a superposition of states (left). The qubit collapses to one state when this is measured (right). The state the qubit will collapse is not known a priori and depends on the basis (or, ``axis'') used in the measurement (dashed arrows). The qubit in the figure is measured in the Z-basis, also known as the computational basis. %\textcolor{violet}{Since the measurement is in z axis, shouldn't it collapse to this axis (either z or -z) in the figure?}
} 
\label{fig:qubit_state}
\end{figure}

%%%%%%%%%%%%%%%%%%%%%
\section{A Brief Introduction to Quantum Networking}
\label{sec:brief_introduction}
A quantum network is a collection of interconnected quantum devices that distribute entangled qubits. In short, a qubit can be a particle, such as an electron or photon, with intrinsic quantum mechanical properties also known as \emph{state}. For example, the state can be the spinning direction of an electron or the polarization of a photon.
A fundamental property of qubits is that their state is undetermined before they are measured, as illustrated schematically in Fig.\ \ref{fig:qubit_state}. Also, the state of two or more qubits can be correlated, which is a phenomenon known as \emph{entanglement}.

Sometimes it is easier to explain what qubits and entanglements are---and why they are useful---by making an analogy. In brief, we can think of a qubit as a closed box containing a ball that changes color randomly over time. The ball's color represents the state of the qubit, which is revealed only when the box is opened, i.e., when the qubit state is measured. For example, when opening the box, the ball could be red or blue with equal probability. We say two boxes are entangled if the colors of the balls inside are correlated. For instance, if the balls are correlated in the
sense that they always have opposite colors, knowing the color of one ball allows us to infer the color of the other. %
Such phenomenon is important because it allows two observers to witness the same random event independently of the distance they are separated. That is, one of the boxes could be on Earth and the other in Mars, and yet, the two balls would have opposite colors after opening the boxes. %Einstein referred to such phenomenon as ``spooky action at a distance'' since it allows two random events (i.e., the balls' colors) to ``coordinate'' faster than the speed of light.   

\begin{figure}
\centering
\includegraphics[width=0.85\columnwidth]{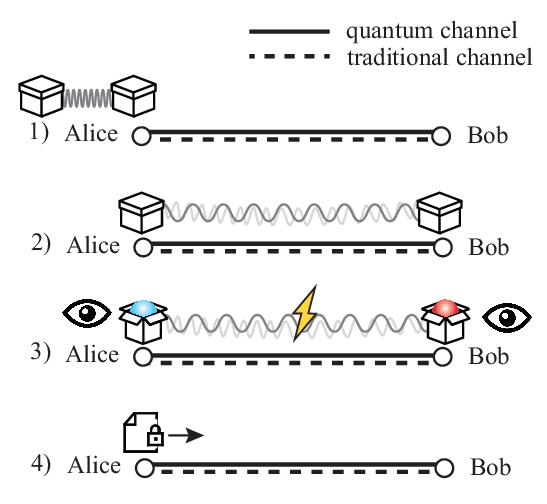}
\caption{Simplified quantum key distribution (QKD) application. Alice and Bob are connected by a quantum link (solid line) and a traditional communication channel (dashed line). 1) Alice generates a pair of entangled qubits. 2)  Alice sends one qubit to Bob. 3) Alice and Bob observe the state of their qubits, and the entanglement breaks. 4) Alices uses the knowledge obtained in step 3 to encrypt a message.}
\label{fig:toy_application}
\end{figure}

We can use the entanglement phenomenon to transmit information securely over a traditional communication network. Consider the simplified QKD protocol shown in Fig.\ \ref{fig:toy_application}, where Alice and Bob are connected by a communication channel (e.g., wireless) and a quantum channel (i.e., fiber optics). First, Alice creates two entangled boxes/qubits and sends one to Bob. Alice and Bob agree to open the boxes and observe the color of their respective balls. Suppose the balls are correlated in the sense that they always have opposite colors. Then, by opening the boxes, Alice and Bob gain access to common knowledge (not information), which they can use for transmitting data (i.e., information) securely over a traditional communication channel. For example, suppose that Alice's ball is blue and that she wants to transmit a $0$ to Bob. Then, she can send the message  $\texttt{\{(\text{blue}: 0), (\text{red}: 1)\}}$ to Bob, which he will interpret as Alice wanting to transmit a $0$---since Bob knows Alice's ball is blue. The communication over the traditional channel is secure because it is not possible for an  eavesdropper that intercepts the message $\texttt{\{(\text{blue}: 0), (\text{red}: 1)\}}$ to infer which bit (i.e., the information) Alice wants to transmit to Bob. 

The toy QKD protocol above requires of course that the quantum channel is secure, which is a requirement that QKD protocols like E91 \cite{E91}  do not have (see also Sec.\ \ref{sec:QKD}). 
Nonetheless, the application brings out three aspects worth emphasizing. First, by sharing qubits, Alice and Bob can observe the \emph{same} random event independently of the distance they are separated. Second, observing the qubits states (i.e., the balls' colors) allows Alice and Bob to acquire \emph{knowledge} that they can  use later for transmitting \emph{information} securely over another channel. And three, no information is transmitted over the quantum channel faster than the speed of light. 
The application also raises the question of how to distribute entangled qubits when clients are not directly connected. However, before we can dive into that (Sec. \ref{sec:qnc}), we need to describe the fundamental properties of qubits and entanglements.

%%%%%%%%%%%%%%%%%%%%%%%%%
\section{Quantum Networking Fundamentals}
\label{sec:fundamentals}

This section describes fundamental concepts of quantum networking that are crucial for distributing entanglement and transmitting quantum information. We begin by explaining the properties of qubits and entanglements using basic probability terminology (Sec.\ \ref{sec:qubitsandentangelments}), and then describe how these enable quantum network control operations for distributing entangled qubits (Sec.\ \ref{sec:quantum_operations}).% Also, we present two quantum applications (QKD and distributed quantum computing) that require of entanglement distribution (Sec.\ \ref{sec:applications}). }}

\subsection{Qubits and quantum entanglement}
\label{sec:qubitsandentangelments}
\begin{figure*}[t!]
\begin{center}
\centering
\includegraphics[width=0.6\textwidth]{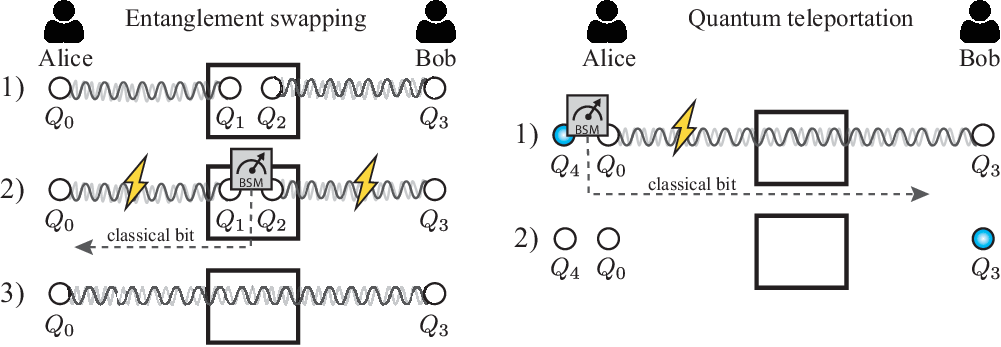}
\caption{Illustrating two quantum network operations. \textbf{Entanglement swapping}: Operation that creates an end-to-end entanglement between two clients (Alice and Bob) via a, for example, quantum repeater (black box). 
The repeater performs a Bell State Measurement (BSM) on qubits $Q_1$ and $Q_2$, which entangles qubits $Q_0$ and $Q_3$ after the exchange of classical information. Fig.\ \ref{fig:circuit} shows the circuit of the entanglement swapping.  \textbf{Quantum teleportation:}  The BSM between $Q_4$ and $Q_0$ transfers the state of qubit $Q_4$ to qubit $Q_3$. %\textbf{Entanglement distillation:}  Operation that combines two or more entanglements to create an entanglement with higher fidelity. 
}
\label{fig:operations}
\end{center}
\end{figure*}

A qubit is the quantum counterpart of a classical bit, but while a classical bit can only exist in a single state (either 0 or 1), a qubit can exist in a superposition or a mix of both states. When a qubit is measured, its state collapses to either 0 or 1 with some probability, where the probability depends on the basis used in the measurement (see example in Fig.\ \ref{fig:qubit_state}). Thus, we can think of a qubit as a random variable
\[
Q:\Omega \to \{0,1\},
\]
where $\Omega$ is the sample space. Measuring a qubit is analogous to drawing a random variable, and we say two or more qubits  are entangled if their outcomes are not independent. Namely, the state in which we observe a qubit gives us information about the state of the others, and vice-versa.
Next, we describe the main properties of qubits and quantum entanglement using a probabilistic lens.

\subsubsection{Qubit} 
\label{sec:qubit_properties}
The probability of observing a qubit in a certain state depends on the basis used in the measurement. Fig.\ \ref{fig:qubit_state} shows, schematically, how an observer sees an electron spinning upwards depending on the ``angle'' of observation (dashed arrows). This property is important for applications such as QKD, as we will show in Sec.\ \ref{sec:QKD}. 
Next, a (qubit) measurement is an active process that changes the probability of observing a qubit in a certain state (i.e., 0 or 1). Or, more specifically, the state of a qubit is random the first time this is observed, and then the qubit remains in the same state in subsequent measurements---if measured in the same basis. This property prevents that the unknown state of a qubit is observed multiple times and, consequently, that a qubit's (unknown) state is copied.

 \subsubsection{Entanglement} 
 \label{sec:entanglement-properties} 
When measuring entangled qubits, the correlation between different measurement outcomes depends on the level of entanglement, formally referred to as entanglement \emph{fidelity}. We say an entanglement has maximum fidelity (or, fidelity one) if the state of a qubit determines the state of the other qubits with probability one. The fidelity of an entanglement may decrease over time due to a phenomenon known as \emph{decoherence} \cite{Bac20}. Decoherence can be caused by noise, which may be the result of different factors, including the distance traveled by the qubits, cross-talk \cite{Ohkura21}, the qubit technology (e.g., superconducting, trapped ions, photons), among others. Finally, an entanglement breaks when one of the qubits is measured, meaning that the qubits' states become independent.

\subsection{Quantum network operations}

Next, we outline  four fundamental operations required for the functioning of quantum networks.

\label{sec:quantum_operations}

\subsubsection{Measurement} This operation was  already mentioned in Sec.\  \ref{sec:qubitsandentangelments} and consists of observing the state of one or multiple qubits with respect to a basis or axis of reference. For example, if an electron is spinning upwards or downwards, which is then mapped to 0 or 1 (see Fig. \ref{fig:qubit_state}). 

\subsubsection{Entanglement swapping} 
\label{sec:swapping} This operation entangles two qubits $Q_0$ and $Q_3$ by using two auxiliary qubits $ Q_1$ and $ Q_2$ that are already entangled with $Q_0$ and $Q_3$  respectively (see Fig.\ \ref{fig:operations}). Importantly, the entanglements between $Q_0$/$Q_2$ and $Q_1$/$Q_3$ break after the swapping operation. This operation is  useful for creating long-distance or end-to-end entanglements with clients that are connected over two or more hops. Also, the operation can be extended to entangle multiple qubits with Greenberger–Horne–Zeilinger (GHZ) states. Fig.\ \ref{fig:circuit} shows an example of a quantum circuit that implements such an operation.

\subsubsection{Teleportation \cite{BBC+93}}
\label{sec:teleportation}
 This operation ``transfers'' the state of a qubit (i.e., its probability distribution) to another qubit. The process is illustrated schematically in Fig.\ \ref{fig:operations}, where the state of the qubit $Q_4$ in Step 1 is transferred to qubit $Q_3$.  This operation is called teleportation because the state of qubit $Q_4$ \emph{before} the operation is equal to the state of $Q_3$ \emph{after} the operation. There are three points worth noting. First, the entanglement between $Q_0$ and $Q_3$ breaks after the teleportation operation is carried out. Second, the state of qubit $Q_4$ \emph{after} the teleportation is not equal to the state of $Q_3$, i.e., the operation does not ``copy'' the state of a qubit. And third, quantum teleportation requires classical communication, and so it is not possible to transmit quantum information faster than the speed of light.

\subsubsection{Entanglement distillation} 
\label{sec:distillation}
This operation combines two or more entanglements to create an entanglement with higher fidelity.
%(see Fig.\ \ref{fig:operations}). 
Distillation plays a crucial role in distributing high-quality entanglements, as the entanglement obtained by sharing a pair of entangled photons over an optical fiber may have low fidelity and/or short decoherence times.  By combining multiple entanglements, it is possible to obtain entanglements with higher fidelity.

\begin{figure*}[ht!]
\centering
\includegraphics[width=\textwidth]{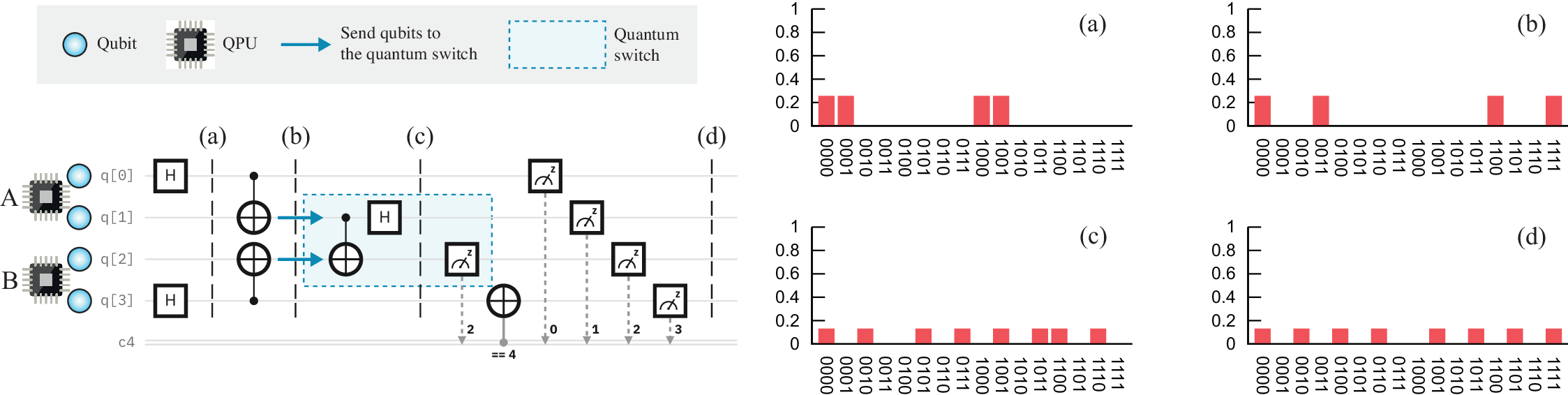}
\caption{Example of a (distributed) quantum circuit.  The circuit input are four qubits in state $0000$ (\texttt{q[3]}, \texttt{q[2]}, \texttt{q[1]}, \texttt{q[0]}), and the goal is to entangle qubit 0 and qubit 3.  The bar graphs on the right-hand-side of the circuit are the qubit state probabilities through the different steps. The first step (a) applies a Hadamard gate to qubits 0 and 3, which makes them be in a superposition of states $0$ and $1$, with probability $1/2$ each. The second step (b) applies a CNOT gate which entangles qubits 0 and 1, and qubits 2 and 3. Observe from the figure that the qubits are indeed entangled since by knowing the state of qubit 0/2, we can infer the state of qubit 1/3. The third step (c) is a Bell State Measurement (CNOT + Hadamard)  which entangles qubits 0 and 2. Finally, the fourth step (d) measures the state of qubit 2, and if this is observed in state 1, it applies a NOT gate to qubit 3. The result is that qubits 0 and 3 are always observed in the same state---see (d). The readers can easily implement such a circuit with IBM's quantum composer \cite{quantum-composer}.}
\label{fig:circuit}
\end{figure*}

\section{Quantum Networking Applications}
\label{sec:applications}

We present two quantum networking applications that require distributing entanglement across nodes in a network (i.e., between Alice and Bob in QKD, and across Quantum Processing Units (QPUs) in distributed quantum computing).

\subsection{Quantum Key Distribution (QKD-E91 \cite{E91})}
\label{sec:QKD}

Consider the setting where Alice is connected to Bob via a quantum channel and another traditional channel (e.g., Fig.\ \ref{fig:toy_application}). Suppose also, for simplicity,  that there is no decoherence and that entanglements have maximum fidelity. That is, by observing the state of one  qubit, we know the states of the other qubits with probability one. The goal is to exchange a \emph{private} sequence of bits, which will be used later as an encryption-key to transfer information securely over a traditional communication channel (e.g., with SSL/TLS). 
The QKD protocol is as follows:
 
\begin{itemize}
\item [(i)] Alice generates $n$ pairs  $(Q_1,  Q_1'), \dots, (Q_n,  Q_n')$ of entangled qubits (i.e., random variables) and sends $( Q_1',\dots,  Q_n')$  to Bob over the quantum channel.  
\item [(ii)] Alice and Bob measure the states of their qubits by selecting bases uniformly at random from a fixed set of options. Recall from Sec.\ \ref{sec:qubit_properties} that the probability of observing a qubit in a certain state depends on the basis used in the measurement.

\item [(iii)] Alice and Bob exchange, over the traditional communication channel, a \emph{sub-sequence} of the qubit measurements and the bases used to measure those. For example, if $n=1000$, Alice and Bob can exchange the sub-sequences of measurements of qubits $(Q_{100},\dots,Q_{199})$ and $( Q_{100}',\dots, Q_{199}')$ and the corresponding bases.
\end{itemize}

With the information in the last step, Alice and Bob can determine if the quantum channel is secure by \emph{comparing the outcomes in the sub-sequence that have been measured with the same bases.} If \emph{all} the outcomes measured with the same bases coincide, Alice and Bob can be confident that the communication over the quantum channel was secure, and use \emph{all} the elements in the sequence measured with the same bases as a private key.\footnote{Their confidence depends on the length of the sub-sequence exchanged.} 
However, if one or more outcomes in the sub-sequence measured with the same bases do not coincide, Alice and Bob can suspect that something went wrong. For instance, an eavesdropper could have intercepted and measured one or more qubits with bases different from the bases used by Alice and Bob. Another possibility could be that the entanglements had low fidelity and that the measurement outcomes did not coincide even though Alice and Bob measured the qubits with the same bases. In either case, the protocol must be restarted. 

\subsection{Distributed quantum computing} 

Increasing the number of qubits in a quantum computer is essential for running programs, such as Shor's algorithm \cite{Sho99} for factoring integers as the product of two prime numbers. However, having a quantum computer with a large number of qubits is challenging as  more qubits often leads to higher levels of noise, which negatively impacts the quality of the qubits' measurements in the circuit's output. 

One of the most promising approaches to increase the number of qubits available and their quality is to interconnect multiple quantum processors, as shown in IBM's quantum roadmap \cite{ibmroadmap}. 
However, such an approach requires not only designing circuits/programs that can be executed in a distributed manner, but also network protocols that enable quantum processors to entangle their qubits (e.g., via entanglement swapping) for quantum information exchange (e.g., through quantum teleportation).

A simple example that illustrates the two points above is the circuit/program in Fig.\ \ref{fig:circuit}. In brief, the circuit correspond to the entanglement swapping protocol described in Sec.\ \ref{sec:swapping}, which can be regarded as a ``distributed'' implementation of a circuit that entangles two qubits (Hadamard + CNOT gates). Processors A and B can store two qubits each, and the goal is to entangle qubit 0 (Processor A) with qubit 3 (Processor B). The entanglement is achieved by entangling, respectively, qubits 0 and 3 with qubits 1 and 2, and, subsequently, transmitting qubits 1 and 2 to a central point (e.g., quantum switch) for a Bell State Measurement (BSM). Observe from the figure that qubits 0 and 3 are entangled at the end of the circuit despite they have never interacted directly. 
The same process can be applied to general quantum programs, typically expressed as a sequence of unitary operators. The goal is to implement these unitary operators in a distributed manner, considering the qubit capacity and communication capabilities of each processor.  
%

%%%%%%%%%%%%%
\section{Quantum Network Control \\ \& Open Research Problems}
\label{sec:qnc}

Distributing entanglement is crucial for realizing the applications described in Sec.\ \ref{sec:applications}. In this section, we present a basic quantum network control model for distributing entanglement independent of the limitations of (near-term) hardware (Sec.\ \ref{sec:network_model}), and then discuss open problems in quantum network control  (Sec.\ \ref{sec:open_problems}). %\textcolor{violet}{This section can be motivated by the previous section by saying why distributing entangleemnt is important (this was motivated before but is a good bridge?)}

\begin{figure*}
\begin{center}
\includegraphics[width=0.5\textwidth]{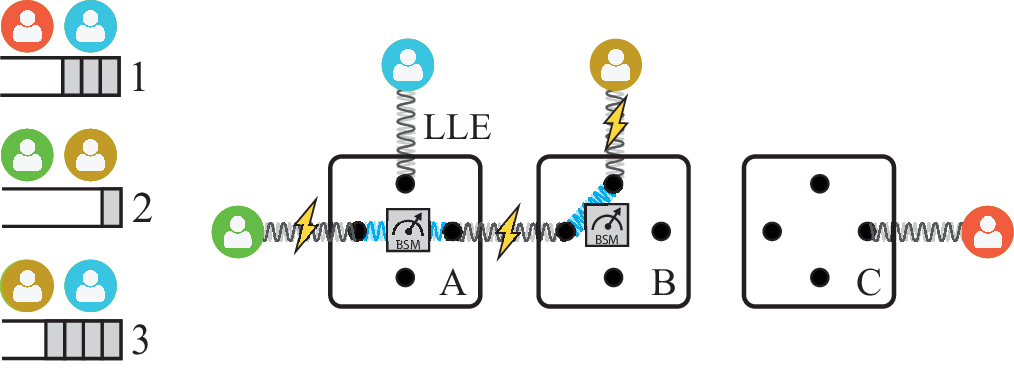}
\caption{Example of a multi-hop quantum network with four clients, three quantum switches (black boxes), and three types of requests. The gray and blue wiggly lines represent LLEs and BSMs, respectively. The requests are stored in separate queues. The figure shows that to serve a request from the second queue, the network controller needs to perform two BSMs with three LLEs.}
\label{fig:qnet}
\end{center}
\end{figure*}

\subsection{Basic quantum network control model}
\label{sec:network_model}

We can model the operation of a quantum network as a discrete time system. In brief, time is divided into slots of equal duration, where in each time slot, the quantum clients (i.e., the nodes) exchange entangled qubits with their neighbors to create \emph{link-level entanglements} (LLEs). LLEs can last for one or multiple time slots depending on the decoherence time of the qubits. In addition, quantum applications (e.g., QKD) generate requests for end-to-end entanglements between clients in each time slot, which are stored in separate queues. 
% The queues in a time slot $t$ evolve as follows
% \begin{align}
% X(t+1) = \left[ X(t) - S(t) \right]^+ + A(t),
% \label{eq:queue_recursion}
% \end{align}
% where $A(t)$ and $S(t)$ are the requests arrivals and service at time $t$ respectively. 
The task of the quantum network controller is to perform entanglement swappings (Sec.\ \ref{sec:swapping}) with the available LLEs to create end-to-end entanglements, which are used to serve the requests in the queues. For example, connecting Alice and Bob as shown in Fig.\ \ref{fig:operations} corresponds to generating an end-to-end entanglement, and serving a request corresponds to using the entanglement to, for example, teleport a qubit (see Sec.~\ref{sec:teleportation} and Fig.\ \ref{fig:operations}). Importantly, the requests that the network controller can serve depends on the available LLEs.

Fig.\ \ref{fig:qnet} shows an example of the connectivity of a multi-hop quantum network in a time slot (with five clients, three quantum switches, and three types of requests). Observe from the figure that the LLEs available (gray wiggly lines) %\textcolor{violet}{Now its the blue wiggly lines, right?} 
constrain the BSMs that the switches can perform, and therefore, the requests that it is possible to serve in a time slot. In particular, it is possible to serve a request from queues 2 and 3 
%\textcolor{violet}{I dont think that it is possible to serve a request in queue 3 since the blue client does not have LLE. Right?} 
but not from queue 1 since switches B and C do not share an LLE. Also, note that the network controller has to decide which of the requests  from queues 2 and 3 to serve since an LLE can only be used to serve one request at a time. Furthermore, unlike traditional communication networks, the links in quantum networks (i.e., LLEs) are consumed to serve requests (see also Sec.\ \ref{sec:swapping}).

We can operate a quantum network with techniques from wired and wireless networks, but we must accommodate the unique characteristics that are not present in traditional networks. In particular, that LLEs act as ``virtual'' links that are consumed to serve requests, and that the quality of those links may degrade rapidly. Furthermore, quantum devices have often limited memory, and there are operations such as entanglement distillation (Sec.~\ref{sec:distillation}) that are unique to quantum networks. In the following section, we discuss some of the challenges and open research questions in quantum network control.

%%%%%%%%%%%
\subsection{Open research problems}
\label{sec:open_problems}

\subsubsection{Characterizing the capacity of a quantum network}
\label{sec:ext_decoherence}

A fundamental research problem is to characterize the maximum number of requests that a quantum network can support, which is often referred to as the network capacity region \cite{9951209}. To address this, we need to consider three important features unique to quantum networks: (i) decoherence, (ii) entanglement distillation, and (iii) quantum memory availability. 

Understanding the network capacity region in relation to these features is crucial for two reasons. Firstly, it will allow us to quantify their impact on the maximum load that the network can support, which is essential for designing quantum networks in a cost-effective manner. Secondly, it will allows us to develop algorithms for operating quantum networks and compare their efficiency against the theoretically achievable maximum performance.

Characterizing the capacity region requires enumerating all the possible ways in which a controller can serve requests. This task is challenging due to the random nature of LLEs and because their ``quality'' is affected by decoherence and entanglement distillation. To address this problem, a first step could be to analyze the three aforementioned features separately. For instance, we could study the impact of decoherence on the capacity region without distillation and quantum memory limitations. Similarly, it would be interesting to explore specific cases that are more amenable to analysis than the general scenario. For example, by studying the cases when entanglements decohere  \emph{rapidly} or \emph{never}, we can lower/upper bound, respectively, the minimum/maximum load a quantum network can support as a function of decoherence. Finally, it is worth focusing on topologies like star or line networks as those will be the building blocks of future quantum networks. %\textcolor{violet}{Do you think maybe we should put references here is what people are doing?}

\subsubsection{Design and operation of quantum overlay networks}
A quantum network can use LLEs to create a “virtual” link between any pair of nodes in the network. This enables us to create a quantum overlay network of arbitrary topology \cite{pouryousef2022quantum}. For instance, suppose that Alice and Bob are connected by a line of quantum repeaters that periodically generates entanglements connecting the line ends. That makes Alice and Bob virtually connected directly but also creates a pool of entanglements that are readily available when needed. Namely, entanglements are accumulated over time to be used, for example, by applications that require many end-to-end entanglements in a brief timeframe.

Realizing quantum overlay networks involves two primary challenges. The first one is determining the optimal placement of cryogenic quantum memories, which can store high-quality entanglements for extended durations. However, these memories are scarce as they need to operate  under low-temperature conditions.
The second challenge is designing entanglement distribution policies that replenish the cryogenic memories and refresh the quality of existing entanglements through distillation techniques. To address these challenges effectively, it is crucial to first characterize the capacity region of a quantum network: we need to quantify how the placement of cryogenic memories will affect the ability of a quantum network to distribute entanglements.

\subsubsection{Entanglement distribution models without queues}

Emerging entanglement distribution approaches are based on queue scheduling techniques that treat requests for end-to-end entanglements as packets in traditional communication networks (e.g., \cite{9951209}).
In particular, requests are accumulated in queues with unlimited storage capacity (see Fig.\ \ref{fig:qnet}), and the goal is to design a policy that keeps the queues stable \cite{GNT06}. While these approaches are well-suited when requests for end-to-end entanglements are classical information (e.g., QKD), they are not appropriate when requests are \emph{information qubits} that have to be teleported.  For example, consider the scenario where qubits at the output of a quantum circuit must be fed into a circuit located in another quantum processor.  Storing \emph{information qubits} for an extended period is generally not possible due to the short decoherence times of superconducting qubits of current quantum computers, but also because of limited quantum memory availability. 

An alternative to queuing models for teleporting qubits could involve employing scheduling techniques inspired by high-speed optical switches, in which photons seamlessly ``flow'' through the network. However, this approach presents two challenges. Firstly, optical switches are typically modeled as bipartite graphs, allowing efficient decomposition of a traffic demand matrix using well-established Birkhoff-based algorithms (e.g., \cite{LML+15}). Unfortunately, these techniques cannot be directly applied to quantum networks, as they do not generally conform to a bipartite graph structure. Nevertheless, we can draw inspiration from these methods. The second challenge is that existing techniques in optical switching assume static network connectivity. However, the link connectivity of a quantum network is dynamic due to the random nature of LLEs, and so existing algorithms must be adapted to account for this distinctive feature.

\section{Concluding Remarks}

Quantum network control is crucial for realizing quantum networking applications such as QKD and distributed quantum computing. However, the field of quantum networking is still in its infancy, despite the progress made in recent years. Commercial quantum networking hardware is not yet available, and it is unclear what capabilities near-term hardware will have---e.g., in terms of coherence times and memory.

In this paper, we have presented an introduction to quantum network control that is independent of a specific technology. However, to fully leverage the potential of quantum networks, we will need to tailor algorithms to the capabilities of the hardware and understand what is theoretically achievable. Among the various challenges in the area, characterizing the capacity of a quantum network depending on the system characteristics (e.g., coherence times, distillation, memory) stands out as a crucial open research problem.

\section{Acknowledgements}
This work was supported by the European Union’s Horizon 2020 Research and Innovation Program under the Marie Skłodowska-Curie under Agreement 795244. 
The research work was supported by the Army Research Office MURI under the project number W911NF2110325 and by the National Science Foundation under project numbers EEC-1941583 CQN ERC and CNS 1955744.

\bibliographystyle{IEEEtran}
\bibliography{references}

\end{document}